\documentclass[12pt]{article}

\usepackage{amstext    }
\usepackage{amsthm    }
\usepackage{a4}
\usepackage[mathscr]{eucal}

\usepackage{amsmath}
\usepackage{amssymb}
\usepackage{amscd}

\numberwithin{equation}{section}

\newtheorem{theorem}{Theorem}[section]

\newtheorem{proposition}[theorem]{Proposition}
\newtheorem{corollary}[theorem]{Corollary}
\newtheorem{lemma}[theorem]{Lemma}
\newtheorem{remark}[theorem]{Remark}

\newcommand{\cali}[1]{\mathscr{#1}}

\newcommand{\dist}{{\rm \ \! dist}}

\newcommand{\ddc}{dd^c}

\newcommand{\Cc}{\cali{C}}
\newcommand{\Dc}{\cali{D}}

\newcommand{\Lc}{\cali{L}}

\newcommand{\Uc}{\cali{U}}

\newcommand{\C}{\mathbb{C}}

\newcommand{\R}{\mathbb{R}}
\renewcommand{\P}{\mathbb{P}}

\title{Exponential mixing for automorphisms on compact K{\"a}hler manifolds}
\author{Tien-Cuong Dinh and Nessim Sibony}

\begin{document}
\maketitle

\begin{abstract}
Let $f$ be a holomorphic automorphism of positive entropy
on a compact K{\"a}hler surface.
We show that the equilibrium measure of $f$ is exponentially mixing.
The proof uses some recent development on the pluripotential
theory. The result also holds 
for automorphisms on compact K{\"a}hler manifolds of higher dimension under a natural
condition on their dynamical degrees. 
\end{abstract}




\section{Introduction}

Let $(X,\omega)$ be a compact K{\"a}hler manifold of dimension $k$ and
$f$ a holomorphic automorphism of $X$. 
The {\it dynamical degree of order $q$} of $f$ is the spectral
radius of the pull-back operator $f^*$ acting on the Hodge cohomology group
$H^{q,q}(X,\C)$. It is denoted by $d_q(f)$ or simply by $d_q$ if
there is no confusion. We have $d_0=d_k=1$ and
if $f^n:=f\circ\cdots\circ f$ ($n$ times)
is the iterate of order $n$ of $f$, then  $d_q(f^n)=d_q^n$.

A theorem by Khovanskii \cite{Khovanskii}, Teissier \cite{Teissier}
and Gromov \cite{Gromov} implies that the sequence $q\mapsto \log d_q$
is concave. So, there are integers $0\leq
p\leq p'\leq k$ such that 
$$1=d_0<\cdots <d_p=\cdots=d_{p'}>\cdots>d_k=1.$$
An instructive example with $p\not=p'$ is a map $f$ on a product
$X=Y\times Z$ of compact K{\"a}hler manifolds such that
$f(y,z)=(g(y),z)$ for $(y,z)\in Y\times Z$. More
interesting examples of maps  preserving a fibration were considered in \cite{Dinh}.

Most dynamical studies on automorphisms of compact
K{\"a}hler manifolds are concentrated on the case where the consecutive dynamical
degrees are distinct, i.e. $p=p'$. Somehow, this condition insures
that the considered dynamical systems have no trivial direction. From now on, we
also assume that $f$ satisfies this natural condition. In  \cite{DinhSibony4, DinhSibony12},
we constructed for $f$ canonical invariant currents (Green currents)
and ergodic invariant
probability measures using the theory of intersection of currents, see
also \cite{Guedj}.

When the operator $f^*$, acting on $H^{p,p}(X,\C)$, admits only one eigenvalue of
maximal modulus, there is 
only one invariant probability measure obtained as the intersection of
a Green $(p,p)$-current of $f$ and a Green $(k-p,k-p)$-current of
$f^{-1}$. We call it the {\it
  equilibrium measure} of $f$. The above eigenvalue
is necessarily equal to $d_p$ and the obtained measure is shown to be 
mixing, hyperbolic and of maximal entropy. The reader finds in \cite{DinhSibony12} and in Section
\ref{section_mixing} below some details. Here is our main theorem.

\begin{theorem} \label{th_main}
Let $f$ be a holomorphic automorphism on a compact K{\"a}hler manifold
$(X,\omega)$ and $d_q$ its dynamical degrees. Assume that there is a degree
$d_p$ strictly larger than the other ones and that
 $f^*$, acting on $H^{p,p}(X,\C)$, admits only one eigenvalue
of maximal modulus $d_p$. 
Then the equilibrium probability measure $\mu$ of
 $f$ is exponentially mixing. More precisely, if  
$\delta$ is a constant such that $\max(d_{p-1},d_{p+1})<\delta<d_p$ and
all the eigenvalues of $f^*$, acting on $H^{p,p}(X,\C)$, except $d_p$,
 are strictly smaller than $\delta$, then
$$|\langle \mu, (\varphi\circ f^n)\psi\rangle
-\langle\mu,\varphi\rangle\langle\mu,\psi\rangle|\leq
A\|\varphi\|_{\Cc^\beta}\|\psi\|_{\Cc^{\beta'}}
(d_p/\delta)^{-n\beta\beta'/8},$$
for all $\Cc^\beta$ function $\varphi$ and all $\Cc^{\beta'}$ function
$\psi$  on $X$ with $0\leq \beta,\beta'\leq 2$. Here, 
$A=A(\beta,\beta',\delta)$ is a constant independent of $\varphi$, $\psi$ and of the
integer $n\geq 0$. 
\end{theorem}

Mixing is equivalent to the property that the left hand
side of the above inequality converges to 0 when $n$ goes to
infinity. 
In the proof of Theorem \ref{th_main}, we use in particular dynamical properties of the map
$F:= (f^{-1},f)$ acting on $X\times X$. A Green $(k,k)$-current of $F$
can be obtained as the limit of $d_p^{-2n}(F^n)^*[\Delta]$, 
where $[\Delta]$ is the current of integration on the
diagonal $\Delta$ of $X\times X$. The speed of convergence 
is the key point in the proof of our result, see Proposition
\ref{prop_green_current} below. The idea was already introduced in
\cite{Dinh1, DinhNguyenSibony, DinhSibony12}. However, the use of the pseudoconvexity of $\C^k$ is no longer valid in the 
compact setting. We will replace it with the use of the H\"older continuity of Green super-potentials.  

Theorem \ref{th_main} still holds under weaker hypothesis: 
all the eigenvalues of maximal modulus of $f^*$, acting on $H^{p,p}(X,\C)$,
 are equal to $d_p$ and the spectral radius $d_p$ of this operator 
 is of multiplicity 1, i.e. $\|(f^n)^*\|\sim d_p^n$. 
The last property can be seen in 
the Jordan form of the
square matrix associated to $f^*$: the Jordan blocks
whose diagonal entries have
modulus $d_p$, are of size $1\times 1$. Very likely, the condition 
$\|(f^n)^*\|\sim d_p^n$ is 
necessary because it insures that the cohomology classes associated
to $d_p^{-2n}(F^n)^*[\Delta]$ converge exponentially fast. Otherwise, we cannot have a good speed of 
convergence for the currents $d_p^{-2n}(F^n)^*[\Delta]$.
In the considered case, the
construction in \cite{DinhSibony12} gives a finite family of
invariant probability measures of
maximal entropy. They are all exponentially mixing. 
 
Consider now an automorphism $f$
of positive entropy on a
compact K{\"a}hler surface $X$. 
Results by Gromov \cite{Gromov1} and Yomdin \cite{Yomdin} say that
the (topological) entropy of $f$ is equal to $\log d_1$. So, $d_1>1$
and the consecutive
dynamical degrees of $f$ are distinct. These automorphisms were
studied by Cantat in \cite{Cantat}. He showed in particular that
all the eigenvalues of $f^*$, acting on
$H^{1,1}(X,\C)$, have modulus 1 except two eigenvalues $d_1$ and
$1/d_1$. So, we can apply Theorem \ref{th_main} and
deduce the following result.

\begin{corollary}
Let $f$ be a holomorphic automorphism of positive entropy on a compact
K{\"a}hler surface $X$. Then the equilibrium measure of $f$ is
exponentially mixing.
\end{corollary}

Note that exponential mixing for polynomial automorphisms was
proved by Nguyen and the authors 
in \cite{Dinh1, DinhNguyenSibony}. We refer to 
Bedford-Kim \cite{BedfordKim}, Keum-Kondo \cite{KeumKondo}, McMullen
\cite{McMullen} and Oguiso \cite{Oguiso}
for interesting examples of automorphisms on compact K{\"a}hler manifolds.

\bigskip\noindent
{\bf Acknowledgement.} The first author wishes to express his gratitude to the Max-Planck Institut f\"ur Mathematik in Bonn
for its hospitality during the preparation of this paper.


\section{Super-potentials of currents} \label{section_sp}

Super-potentials were introduced by the authors 
in order to develop a calculus on positive closed
currents. We recall some basic properties and refer
to \cite{DinhSibony12} for details.

Let $\Dc_p$ denote the real space generated by
positive closed $(p,p)$-currents on $X$. If $S$ is a current in
$\Dc_p$, define the {\it norm} $\|S\|_\ast$ of $S$ by
$$\|S\|_\ast:=\min \|S^+\|+\|S^-\|$$
where the minimum is taken over the positive closed currents $S^\pm$
with $S=S^+-S^-$. Here, $\|S^\pm\|$ denote the {\it mass} of $S^\pm$ which
are defined by
$$\|S^\pm\|:=\langle S,\omega^{k-p}\rangle.$$ 
Observe $\|S^\pm\|$ depend only 
on the cohomology classes of $S^\pm$ in $H^{p,p}(X,\R)$. We say that a
subset of $\Dc_p$ is {\it $\ast$-bounded} if it is bounded for the
$\|\ \|_\ast$-norm. 
Let $\Dc_p^0$ denote the subspace of currents $S$ in $\Dc_p$ whose
classes $\{S\}$ in $H^{p,p}(X,\R)$ are zero. 

We consider on
$\Dc_p$ and $\Dc_p^0$ the following {\it topology}: a sequence $(S_n)$
in $\Dc_p$ or $\Dc_p^0$ converges to a current $S$ if $S_n$ converge
to $S$ in the sense of currents
and if $\|S_n\|_\ast$ are bounded by a constant independent of
$n$. Smooth forms are dense in $\Dc_p$ and $\Dc_p^0$ for this topology.

For any $0<l<\infty$, we can associate to $\Dc_p$ a {\it norm} $\|\
\|_{\Cc^{-l}}$ and a {\it distance} $\dist_l$ defined by
$$\|S\|_{\Cc^{-l}}:=\sup_{\|\Phi\|_{\Cc^l}\leq 1}|\langle S,\Phi\rangle|\quad
\mbox{and}\quad \dist_l(S,S'):=\|S-S'\|_{\Cc^{-l}},$$
where $\Phi$ is a smooth test form of bidegree $(k-p,k-p)$ on $X$. The
weak topology on each $\ast$-bounded subset of $\Dc_p$ coincides with
the topology induced by $\|\ \|_{\Cc^{-l}}$. If $0<l<l'<\infty$ are
two constants, then on each $\ast$-bounded subset of $\Dc_p$ we have
$$\dist_{l'}\leq \dist_l\leq c_{l,l'}(\dist_{l'})^{l/l'}$$
for some positive constant $c_{l,l'}$.

The super-potential of a current $S$ in $\Dc_p$ is a canonical linear
function defined, under some normalization, on the smooth forms
in $\Dc_{k-p+1}^0$. It plays the same role as the potentials of
positive closed $(1,1)$-currents which are quasi-p.s.h. functions. 

Let $\alpha=(\alpha_1,\ldots,\alpha_h)$ with $h:=\dim H^{p,p}(X,\R)$
be a fixed family of real smooth
closed $(p,p)$-forms such that the family of classes
$\{\alpha\}=(\{\alpha_1\},\ldots,\{\alpha_h\})$ is a basis of
$H^{p,p}(X,\R)$. Let $R$ be a current in $\Dc^0_{k-p+1}$. Since the cohomology class of
$R$ is zero, there is 
a real $(k-p,k-p)$-current $U_R$ such that 
$\ddc U_R=R$. We call $U_R$ a {\it potential of $R$}. 
Adding to $U_R$ a suitable closed form allows to assume
that $\langle U_R,\alpha_i\rangle=0$ for
$i=1,\ldots,h$ and we say that $U_R$ is {\it $\alpha$-normalized}. When
$R$ is smooth, we can choose $U_R$ smooth and the {\it $\alpha$-normalized
super-potential} $\Uc_S$ of $S$ is defined by
$$\Uc_S(R):=\langle S, U_R\rangle.$$
The definition does not depend on the choice of $U_R$.

When the
function $\Uc_S$ extends continuously to $\Dc_{k-p+1}^0$ for the considered topology,
we say that $S$ has a {\it continuous} super-potential.
If $S$ is in $\Dc_p^0$ then $\Uc_S$ does not depend on the choice of $\alpha$; if
moreover $S$ is smooth, it has a continuous super-potential and we
have the formula
$$\Uc_S(R)=\Uc_R(S),$$
where $\Uc_R$ is the super-potential of $R$ which is also independent
of the normalization. We can extend the above equality to the case
where $S$ has a continuous super-potential.

We say that $\Uc_S$ is {\it
  $(l,\lambda,M)$-H{\"o}lder continuous} if it is continuous and if 
$$|\Uc_S(R)|\leq M\|R\|_{\Cc^{-l}}^\lambda$$ 
for $R\in\Dc_{k-p+1}^0$ with $\|R\|_\ast\leq 1$,
where $l>0$, $0<\lambda\leq 1$ and $M\geq 0$ are constants.
If $l'>0$ is another constant, the above comparison between $\dist_l$
and $\dist_{l'}$ implies that when $\Uc_S$ is $(l,\lambda,M)$-H{\"o}lder
continuous, it is $(l',\lambda',M')$-H{\"o}lder continuous for some
constants $\lambda'$ and $M'$ which are independent of $S$.

Here is the main result in this section. It improves Theorem 3.2.6 in
\cite{DinhSibony12} and can be seen as a version of the classical exponential
estimates for p.s.h. functions.

\begin{proposition} \label{prop_exp_estimate}
Let $R$ be a current in $\Dc_{k-p+1}^0$ with
$\|R\|_\ast\leq 1$ such that its
super-potential $\Uc_R$ is $(2,\lambda,M)$-H{\"o}lder continuous. 
Then there is a constant $A>0$ independent of $R,\lambda$ and $M$ such
that the super-potential $\Uc_S$ of $S$ satisfies
$$|\Uc_S(R)|\leq A (1+\lambda^{-1}\log^+M),$$
for any current $S$ in $\Dc_p^0$ with $\|S\|_\ast\leq 1$, where $\log^+:=\max(0,\log)$. 
\end{proposition}

We will use a family of linear regularizing operators
$\Lc_\theta:\Dc_p^0\rightarrow\Dc_p^0$  introduced in
\cite{DinhSibony12} with $\theta$ in $\P^1=\C\cup\{\infty\}$. 
Let us recall some properties of $\Lc_\theta$. Fix a constant $c>0$ 
large enough which depends only on the geometry of $(X,\omega)$.

The operators
$\Lc_\theta$ are continuous for the considered
topology on $\Dc_p^0$
and are bounded for the $\|\ \|_\ast$-norm,
i.e. $\|\Lc_\theta(S)\|_\ast\leq c\|S\|_\ast$ with $c>0$ independent of $\theta$ and $S$.
We have $\Lc_0(S)=S$ and 
$$\dist_2(S,\Lc_\theta(S))\leq c\|S\|_\ast|\theta|.$$
Moreover, $\Lc_\theta=\Lc_\infty$ for $|\theta|\geq 1$.

Let $p\geq 1$ be a constant and let $q\geq 1$
such that when $p<k+1$, $1/q=1/p-1+k/(k+1)$ and $q=\infty$ when $p\geq
k+1$. We always have $\|\Lc_\infty(S)\|_{L^1}\leq c\|S\|_\ast$.
If $S$ is an $L^p$ form, $p\geq 1$, then $\Lc_\infty(S)$ is an $L^q$
form satisfying
$$\|\Lc_\infty(S)\|_{L^q}\leq c \|S\|_{L^p}.$$
Here, $c>0$ is a constant large enough.
Recall also that the function $u_S(\theta):=\Uc_{\Lc_\theta(S)}(R)$ 
is continuous and is constant out of the unit disc. It satisfies
$$\|\ddc u_S(\theta)\|\leq c\|S\|_\ast\|R\|_\ast.$$
The last properties hold
for $R$ smooth and extend by continuity to currents $R$
with a continuous super-potential.    

\bigskip

\noindent
{\bf Proof of Proposition \ref{prop_exp_estimate}.}
For $R$ and $S$ as in the proposition, we have $\|S\|_\ast\leq 1$ and
$\|R\|_\ast\leq 1$. Multiplying $S$ by a constant allows to assume
that $\|S\|_\ast\leq c^{-k-3}$.
Define $S_0:=S$ and $S_{i+1}:=\Lc_\infty(S_i)$ for
$0\leq i\leq k+1$. Define also $u_i(\theta):=\Uc_{\Lc_\theta(S_i)}(R)$ and
$m_i:=u_i(0)=u_{i-1}(\infty)$. Using inductively the above estimates, we get
$\|S_i\|_\ast\leq 1/c$, $\|\ddc u_i\|\leq 1$ and
$\|S_{k+2}\|_{L^\infty}\leq 1$. The last inequality implies that $|m_{k+2}|$
is bounded by a constant independent of $S,R$. Indeed, $R$
always admits a potential $U_R$ of bounded
$L^1$-norm and we have $m_{k+2}=\langle S_{k+2},U_R\rangle$.

We need to show that $|m_0|\leq  A (1+\lambda^{-1}\log^+ M)$ for
some constant $A>0$. For this purpose, we can assume that $M>1$ and 
it is enough to check that
$|m_i-m_{i+1}|\leq A (1+\lambda^{-1}\log M)$ for
some constant $A>0$. We have $m_i-m_{i+1}=v_i(0)$ where
$v_i:=u_i-m_{i+1}$. 
The above properties of $u_i$ imply that 
$v_i$ are continuous, vanish outside the unit disc and satisfy
$\|\ddc v_i\|\leq 1$. The classical exponential estimates for
subharmonic functions imply that $\|e^{|v_i|}\|_{L^1(\P^1)}\leq
c$ for some universal constant $c>0$, 
see \cite[Lemma 2.2.4]{DinhSibony12} and \cite[Th. 4.4.5]{Hormander}. We then deduce that there is a
$\theta$ satisfying $|\theta|\leq M^{-1/\lambda}$ and $|v_i(\theta)|\leq (A-1)+A\lambda^{-1}\log M$
for a fixed constant $A$ large enough. Finally, using the H{\"o}lder
continuity of $\Uc_R$, we get
\begin{eqnarray*}
|v_i(0)-v_i(\theta)| & = &  |\Uc_{S_i}(R)-\Uc_{\Lc_\theta(S_i)}(R)|
= |\Uc_R(S_i)-\Uc_R(\Lc_\theta(S_i))|\\
& \leq &  M\dist_2(S_i,\Lc_\theta(S_i))^\lambda
\leq M|\theta|^\lambda \leq 1.
\end{eqnarray*}
Therefore, $|v_i(0)|\leq A (1+\lambda^{-1}\log M)$. This completes the
proof.
\hfill $\square$

\section{Convergence towards Green currents} \label{section_green}

Let $f$, $d_q$ and $\delta$ be as in Theorem \ref{th_main}. 
Fix a constant $\delta_0<\delta$, close enough to $\delta$, 
so that $\delta_0$ satisfies also the same properties as $\delta$.
We recall
some known facts and refer to \cite{DinhSibony12} for details.
By Poincar\'e duality, the dynamical degree $d_q$ of $f$ is equal to
the degree $d_{k-q}(f^{-1})$ of $f^{-1}$. Since the mass of a positive
closed current can be computed cohomologically, if $S$ is in $\Dc_q$
and $R$ is in $\Dc_{k-p+1}$,
we have $\|(f^n)^*(S)\|_\ast\leq cd_p^n\|S\|_\ast$ and 
$\|(f^n)_*(R)\|_\ast\leq c\delta_0^n\|R\|_\ast$ for some constant
$c>0$ independent of $S,R$ and $n$.

By Perron-Frobenius theorem, the eigenspace $H$ associated to the eigenvalue $d_p$ of $f^*$ acting on
$H^{p,p}(X,\R)$ is a real line. Therefore, $d_p^{-n}(f^n)^*$
converge to a linear operator $L_\infty:H^{p,p}(X,\R)\rightarrow
H$. Under the hypothese of Theorem \ref{th_main}, it is easy to deduce
that on $H^{p,p}(X,\R)$
$$\|d_p^{-n}(f^n)^*-L_\infty\|\leq c (d/\delta_0)^{-n}$$ 
for some constant $c>0$. A {\it Green $(p,p)$-current} $T_+$ of $f$ is a
non-zero positive closed $(p,p)$-current invariant under
$d_p^{-1}f^*$, i.e. $f^*(T_+)=d_pT_+$. 
Its cohomology class $\{T_+\}$ generates the real line $H$. 
Moreover, it is known \cite{DinhSibony12} that $T_+$ is the unique positive closed current in $\{T_+\}$.
So, if $S$ is a current in $\Dc_p$, then $d_p^{-n}(f^n)^*(S)$ converge to a
multiple of $T_+$. Here is the main result of this section.

\begin{proposition} \label{prop_green_current}
Let $f,d_q,\delta$ be as in Theorem \ref{th_main} and $S$ a
current in $\Dc_p$. Let $r$ be the constant such that $d_p^{-n}(f^n)^*(S)$
converge to $rT_+$. Let $R$ be a current in $\Dc_{k-p+1}^0$ with
$\|R\|_\ast\leq 1$ whose super-potential $\Uc_R$ is
$(2,\lambda,1)$-H{\"o}lder continuous. Let $\Uc_+$, $\Uc_n$ be the $\alpha$-normalized
super-potentials of $T_+$ and of $d_p^{-n}(f^n)^*(S)$.
Then 
$$|\Uc_n(R)-r\Uc_+(R)|\leq A(d/\delta)^{-n}$$
where $A>0$ is a constant independent of $R$ and of $n$.  
\end{proposition}

We first prove the following lemma.

\begin{lemma} \label{lemma_iterate_form}
Let $R$ be a current in $\Dc_{k-p+1}^0$ whose
  super-potential $\Uc_R$ is $(2,\lambda,M)$-H\"older continuous.  
Then, there is a constant $A_0\geq 1$ independent
of $R,\lambda,M$ such that the super-potential $\Uc_{f_*(R)}$ of $f_*(R)$ is 
$(2,\lambda,A_0M)$-H\"older continuous.
\end{lemma}
\proof
Let $T$ be a current in $\Dc^0_p$ such that $\|T\|_\ast\leq 1$.  
We have seen that $\|f^*(T)\|_\ast\leq
c$ for some constant $c\geq 1$ independent of $T$. Define $T':=c^{-1}
f^*(T)$. 
If $T$ is smooth and $U_T$ is a smooth potential of $T$, then
$f^*(U_T)$ is a smooth potential of $f^*(T)$ and we have 
$$\Uc_{f_*(R)}(T)=\langle f_*(R),U_T\rangle = \langle R,
f^*(U_T)\rangle = \Uc_R(f^*(T)).$$
Since $\Uc_R$ is continuous and smooth forms are dense in $\Dc^0_p$, 
we deduce that $\Uc_{f_*(R)}$ is
continuous and $\Uc_{f_*(R)}(T)=\Uc_R(f^*(T))$ for every $T$ in
$\Dc^0_p$. Therefore,
$$|\Uc_{f_*(R)}(T)|=c|\Uc_R(T')|\leq cM\|T'\|_{\Cc^{-2}}^\lambda.$$

Now, it is enough to show that $\|f^*(T)\|_{\Cc^{-2}}\leq
c'\|T\|_{\Cc^{-2}}$ for some constant $c'>0$. Consider test
$(k-p,k-p)$-forms $\Phi$ such that $\|\Phi\|_{\Cc^2}\leq 1$. Since
$f^{-1}$ is smooth, there is a constant $c'>0$ such that
$\|f_*(\Phi)\|_{\Cc^2}\leq c'$. It follows that
$$\|f^*(T)\|_{\Cc^{-2}}=\sup_\Phi|\langle f^*(T),\Phi\rangle|
=\sup_\Phi |\langle T,f_*(\Phi)\rangle|\leq c'\|T\|_{\Cc^{-2}}.$$
This completes the proof. 
\endproof

\noindent
{\bf Proof of Proposition \ref{prop_green_current}.}
We have seen that
$\|d_p^{-n}(f^n)^*-L_\infty\|\lesssim (d/\delta_0)^{-n}$ on
$H^{p,p}(X,\R)$. So, the computation in \cite[Lemma 4.2.3]{DinhSibony12} shows that
if $\Uc_S$ is continuous, $|\Uc_n(R)-r\Uc_+(R)|\lesssim
(d/\delta)^{-n}$. Therefore, subtracting from $S$ a smooth closed
$(p,p)$-form allows to assume that $\{S\}=0$ and hence $r=0$. 

Define $R_n:=c^{-1}\delta_0^{-n}(f^n)_*(R)$ where $c\geq 1$ is a fixed
constant large enough. We have $\|R_n\|_\ast\leq 1$. 
Lemma  \ref{lemma_iterate_form}
implies by induction that
$\Uc_{R_n}$ is $(2,\lambda,A_0^n)$-H\"older continuous.
As in the proof of this lemma, we obtain
$\Uc_n(R)=c(d_p/\delta_0)^{-n}\Uc_S(R_n)$. 
Finally, we deduce from Proposition \ref{prop_exp_estimate} that
$$| \Uc_n(R)|=c(d_p/\delta_0)^{-n}|\Uc_S(R_n)|\lesssim n
(d_p/\delta_0)^{-n}.$$
The result follows.
\hfill $\square$

\section{Exponential mixing} \label{section_mixing}

In this section, we prove Theorem \ref{th_main}. Theory of
interpolation between the Banach spaces $\Cc^0$ and $\Cc^2$
\cite{Triebel} implies that it is enough to consider the case $\beta=\beta'=2$,
see \cite{Dinh1,DinhNguyenSibony} for details. Assume now that
$\varphi$ and $\psi$ are $\Cc^2$ functions such that
$\|\varphi\|_{\Cc^2}\leq 1$ and $\|\psi\|_{\Cc^2}\leq 1$. Subtracting
from $\psi$ a constant allows to assume also that
$\langle\mu,\psi\rangle=0$. We have to show that
$$|\langle \mu,(\varphi\circ f^n)\psi\rangle|\lesssim (d/\delta)^{-n/2}.$$
We only need to consider the case where $n$ is even. Indeed, if $n$ is
odd, we can replace $\varphi$ with $\varphi\circ f$ and deduce the
result from the first case. So, it is enough to check that 
$$|\langle \mu,(\varphi\circ f^{2n})\psi\rangle|\lesssim (d/\delta)^{-n}.$$

We will apply Proposition \ref{prop_green_current} to the automorphism $F$ of $X\times
X$ defined by $F(x,y):=(f^{-1}(x),f(y))$. By K{\"u}nneth formula
\cite[Th. 11.38]{Voisin}, there is a canonical isomorphism
$$H^{q,q}(X\times X,\C)=\bigoplus_{s+r=q} H^{s,r}(X,\C)\otimes
H^{r,s}(X,\C).$$
It is not difficult to see that $F^*$ preserves the above
decomposition. So, the dynamical degree of order $k$ of $F$ is equal
to $d_p^2$. 
It was shown in \cite{Dinh} that the spectral radius of $f^*$ on
$H^{r,s}(X,\C)$, which is also the spectral radius of $f_*$ on
$H^{k-r,k-s}(X,\C)$, 
is smaller or equal to $\sqrt{d_rd_s}$. Therefore,
the dynamical degrees and the eigenvalues of
$F^*$ on $H^{k,k}(X\times X,\R)$, except $d_p^2$, are strictly smaller than $d_p\delta_0$.  
So, we can apply  Proposition \ref{prop_green_current} to $F$. 

Let $[\Delta]$ denote the positive closed $(k,k)$-current associated
to the diagonal $\Delta$ of $X\times X$. Recall $\mu$ is the 
wedge-product $T_+\wedge T_-$ of a Green $(p,p)$-current $T_+$ associated to $f$ and a 
Green $(k-p,k-p)$-current $T_-$ associated to $f^{-1}$.
We have $f^*(T_+)=d_p T_+$ and $f_*(T_-)=d_p T_-$. Hence,
$F_*(T_+\otimes T_-)=d_p^2 T_+\otimes T_-$. We deduce from the
uniqueness of Green currents that any  Green $(k,k)$-current of $F^{-1}$ is a multiple of $T_+\otimes T_-$. 
In particular, it has a H{\"o}lder
continuous super-potential. 

Recall that $\|\varphi\|_{\Cc^2}\leq 1$ and $\|\psi\|_{\Cc^2}\leq 1$. Define
$\Phi(x,y):=\varphi(x)\psi(y)$. Since the $\Cc^2$-norm of this
function is bounded,
$\ddc\Phi$ is a current in $\Dc^0_2(X\times X)$ with bounded
$\|\ \|_\ast$-norm. 
If $\Uc$ is its super-potential and $T$ is a current in
$\Dc_{2k}^0(X\times X)$, then $\Uc(T)=\langle\Phi,T\rangle$.
Clearly, $\Uc$ is  $(2,1,M)$-H{\"o}lder continuous for some constant $M>0$
independent of $\varphi,\psi$.
By Proposition 3.4.2 in \cite{DinhSibony12}, the wedge-product of 
currents with H\"older continuous super-potentials has also a H\"older
continuous super-potential. We deduce from the proof of that
proposition and the comparison between the distances $\dist_l$ that  $R:=(T_+\otimes
T_-)\wedge \ddc\Phi$ is a current in $\Dc^0_{k+1}(X\times X)$ with
 a $(2,\lambda,M')$-H{\"o}lder continuous
super-potential for some constants $\lambda, M'$ independent of
$\varphi,\psi$. This current $R$ has also a bounded $\|\ \|_\ast$-norm.
Multiplying $\varphi$ by a constant allows us to assume
that $\|R\|_\ast\leq 1$ and $M'=1$. 

Proposition \ref{prop_green_current} applied to $F$, $[\Delta]$ instead of $f$, $S$
yields 
$$|\Uc_{d_p^{-2n}(F^n)^*[\Delta]}(R)-m|\lesssim (d_p/\delta)^{-n}\quad
\mbox{where}\quad  m:=\lim_{n\rightarrow\infty}  \Uc_{d_p^{-2n}(F^n)^*[\Delta]}(R).$$
On the other hand, since $F_*(T_+\otimes T_-)=d_p^2T_+\otimes T_-$ and
since Green currents are well approximated by smooth forms, the
following calculus holds (see
\cite{DinhSibony12})
\begin{eqnarray*}
\Uc_{d_p^{-2n}(F^n)^*[\Delta]}(R) & = & \langle
d_p^{-2n}(F^n)^*[\Delta],\Phi(T_+\otimes T_-)\rangle \\
& = & \langle [\Delta], d_p^{-2n}(\Phi\circ F^{-n})(F^n)_*(T_+\otimes
T_-)\rangle\\
& = & \langle [\Delta], (\Phi\circ F^{-n})T_+\otimes T_-\rangle\\
& = & \langle (T_+\otimes T_-)\wedge [\Delta], \Phi\circ
F^{-n}\rangle.
\end{eqnarray*}
The same arguments and the fact that  $\mu=T_+\wedge T_-$ is invariant yield
$$\Uc_{d_p^{-2n}(F^n)^*[\Delta]}(R) =  \langle T_+\wedge T_-,(\varphi\circ f^n)(\psi\circ
f^{-n})\rangle
=  \langle \mu,(\varphi\circ f^{2n})\psi\rangle.$$
We deduce from the mixing of $\mu$ that the last
integral tends to 0 since $\langle\mu,\psi\rangle=0$. Therefore, we have $m=0$. This together with the
above estimate on the super-potential of $d_p^{-2n}(F^n)^*[\Delta]$ 
implies that 
$$| \langle \mu,(\varphi\circ f^{2n})\psi\rangle|\lesssim
(d_p/\delta)^{-n},$$ 
and completes the proof of Theorem \ref{th_main}.

\begin{remark} \rm
Let $\delta_+\geq d_{p-1}$ (resp. $\delta_-\geq d_{p+1}$) denote the
smallest number such that 
the eigenvalues of $f^*$ acting
on $H^{p,p}(X,\C)$, except $d_p$, are of modulus smaller than or equal
to $\delta_+$ (resp. $\delta_-$). Theorem
\ref{th_main} still holds for any $\delta$ such that
$${2\log\delta_+\log\delta_-\over \log\delta_+
  +\log\delta_-}<\log\delta<\log d_p.$$
Indeed, there are positive integers $l,m$ such that
$$\max\big(\delta_+^l,\delta_-^m\big)<\delta^{l+m\over 2}$$
and it is enough to follow the proof of Theorem \ref{th_main} where we
replace $F$ 
with the automorphism $(f^{-l},f^m)$.
The details are
left to the reader.
\end{remark}


T.-C. Dinh, UPMC Univ Paris 06, UMR 7586, Institut de
Math{\'e}matiques de Jussieu, F-75005 Paris, France. {\tt
  dinh@math.jussieu.fr}, {\tt http://www.math.jussieu.fr/$\sim$dinh}

\

\noindent
N. Sibony,
Universit{\'e} Paris-Sud, Math{\'e}matique - B{\^a}timent 425, 91405
Orsay, France. {\tt nessim.sibony@math.u-psud.fr}

\end{document}